\documentclass[12pt]{amsart}%
\usepackage{amsmath}
\usepackage{graphicx}
\usepackage{amscd}
\usepackage{geometry}
\usepackage{amsfonts}
\usepackage{amssymb}
\usepackage[mathscr]{eucal}%
\theoremstyle{plain}

\numberwithin{equation}{section}

\begin{document}

\title[ ]{Approximating Certain Cell-Like Maps by Homeomorphisms$^\ddag$}

\author{Robert D. Edwards}

\address{Department of Mathematics, UCLA, Box 951555 Los Angeles, CA 90095-1555}

\email{rde@math.ucla.edu}

\thanks{$^\ddag$From photocopied handwritten manuscript labeled \textquotedblleft Preliminary version, 21 July 77, photocopy\textquotedblright}.

\date{}

\begin{abstract}
Given a proper map $f:M\rightarrow Q$, having cell-like point-inverses, from a
manifold-without-boundary $M$ onto an ANR $Q$, it is a much-studied problem to
find when $f$ is approximable by homeomorphisms, i.e., when the decomposition
of $M$ induced by $f$ is shrinkable (in the sense of Bing). If dimension
$M\geq5$, J. W. Cannon's recent work focuses attention on whether $Q$ has the
\emph{disjoint disc property} (which is: Any two maps of a 2-disc into $Q$ can
be homotoped by an arbitrarily small amount to have disjoint images; this is
clearly a necessary condition for $Q$ to be a manifold, in this dimension
range). This paper establishes that such an $f$ is approximable by
homeomorphisms whenever dimension $M\geq5$ and $Q$ has the disjoint disc
property. As a corollary, one obtains that given an arbitrary map
$f:M\rightarrow Q$ as above, the stabilized map $f\times id(\mathbb{R}%
^{2}):M\times\mathbb{R}^{2}\rightarrow Q\times\mathbb{R}^{2}$ is approximable by
homeomorphisms. The proof of the theorem is different from the proofs of the
special cases in the earlier work of myself and Cannon, and it is quite
self-contained. This work provides an alternative proof of L. Siebenmann's
Approximation Theorem, which is the case where $Q$ is given to be a manifold.
\\\\
\\\\
\\\\
\\\\
\textsc{The Edwards' Manuscript Project.}  This article is one of three highly influential articles on the topology of manifolds written by Robert D. Edwards in the 1970's but never published.
This article \textquotedblleft Approximating certain cell-like maps by homeomorphisms\textquotedblright presents the definitive theorem on the recognition of high-dimensional manifolds among resolvable generalized manifolds.  The theorem says that every cell-like map from an $n$-manifold ($n \geq 5$) to an ANR is approximable by homeomorphisms provided that the ANR has the disjoint disc property. (This work garnered Edwards an invitation to give a one-hour plenary address to the 1978 International Congress of Mathematicians.)  The second article \textquotedblleft Suspensions of homology spheres\textquotedblright presents the initial solutions
of the fabled Double Suspension Conjecture. The third article \textquotedblleft Topological
regular neighborhoods\textquotedblright develops a comprehensive theory of regular neighborhoods
of locally flatly embedded topological manifolds in high dimensional
topological manifolds. The manuscripts of these three articles have circulated
privately since their creation. The organizers of the Workshops in Geometric
Topology (http://sites.coloradocollege.edu/geometrictopology2016/) with
the support of the National Science Foundation have facilitated the preparation
of electronic versions of these articles to make them publicly available. Preparation of the
electronic manuscripts was supported by NSF Grant DMS-0407583.  Final editing was carried out by Fredric Ancel, Craig Guilbault and Gerard Venema.
\end{abstract}

\maketitle

\section{Introduction}

\noindent \textbf{Approximation Theorem:} Suppose $f: M \rightarrow Q$ is a proper cell-like map
from a manifold-without-boundary $M$ onto an ANR $Q$, and suppose that $Q$ has
the disjoint disc property and that $\dim M \geq5$. Then $f$ is arbitrarily
closely approximable by homeomorphisms. Stated another way, the decomposition
of $M$ induced by $f$ is shrinkable (in the sense of Bing).

\section{Shrinking Certain 0-dimensional decompositions}

The goal of this section is to prove the following theorem. We are assuming for proofs that M is compact.
\newline%

\noindent \textbf{$\boldsymbol{0}$-Dimensional Shrinking Theorem.} \textit{Suppose $f : M
\rightarrow Q$ is a cell-like map of a manifold $M$ onto a quotient space $Q$,
and suppose }

\begin{itemize}
\item[(1)] \textit{the image of the nondegeneracy set of $f$ has dimension 0,
and }

\item[(2)] \textit{the nondegeneracy set of $f$ has codemension $\geq3$ in
$M$. }
\end{itemize}
\textit{Then the decomposition of $M$ induced by $f$ is shrinkable; that is,
$f$ is arbitrarily closely approximable by homeomorphisms.} 
\newline%
\newline\noindent\textbf{Note.} There is no dimension restriction on $M$ in
this theorem (or anywhere in this section). Nor is there any restriction on
the $closure$ in $Q$ of the image of the nondegeneracy set - it may
be all of $Q$. This is why the theorem will be useful, in sections 3 and 4.
\\\\
\indent The above theorem will be derived from the following theorem, which amounts to
the special case when the decomposition is countable and null.
\newline%
\newline\textbf{Countable Shrinking Theorem.} \textit{Suppose $f : M
\rightarrow Q$ is a cell-like map of a manifold $M$ onto a quotient space $Q$,
and suppose }

\begin{itemize}
\item[(1)] \textit{the nontrivial point-inverses of $f$ comprise a countable
null collection, where \underline{null} means that their diameters tend to 0,
and }

\item[(2)] \textit{each nontrivial point-inverse of $f$ has codimension
$\geq3$ in M. }
\end{itemize}

\textit{Then the decomposition of $M$ induced by $f$ is shrinkable; that is,
$f$ is arbitrarily closely approximable by homeomorphisms.} 
\\\\
\textbf{Notes.} \vspace{6pt}

\indent (1). Both of these theorems are false when codimension $\geq3$ is replaced by
codimension $\geq2$, even assuming $f$ is cellular. In dimension 3, Bing's
countable planar-Knaster-continua decomposition \cite{Bi2} provides a
counterexample to the two theorems. In dimensions $\geq4$, Eaton's generalized
dogbone space \cite{Ea} provides a counterexample to the first theorem, and a
modification of this example, implicit in the first proof below, provides a
counterexample to the second theorem. \vspace{6pt}

\indent (2). As an incidental fact, recall that in either theorem, even without conditions (2),
the quotient space $Q$ is necessarily an ANR, since by hypothesis $Q$ is a
union of two finite dimensional subsets (namely, the image of the
nondegeneracy set of $f$, and its complement), hence $Q$ is finite dimensional.
(See \cite{Hu-Wa}.)
\newline
\newline\textbf{\textit{Proof that Countable Shrinking Theorem $\Rightarrow$ 0-Dimensional Shrinking Theorem.}}  
In brief, the idea is to tube together the nontrivial point inverses of $f$ in
such a manner as to come up with a countable null collection. (This sort of
operation, for the closed-0-dimensional situation, was done in \cite{Ea-Pi} and
Lemma 2 of \cite{Ed-Mi}, and probably elsewhere, too.)

\indent The rationale is this. Given $f : M \rightarrow Q$ as in the first
theorem, suppose that for any $\epsilon> 0$ one can come up with a quotient
map $g : Q \rightarrow Q_{\#}$, with point-inverses of diameter $< \epsilon$,
such that the composition $gf : M \rightarrow Q_{\#}$ is approximable by
homeomorphisms. (In our case, this will be because $gf$ satisfies the Countable
Shrinking Theorem). Then $f$ is approximable by homeomorphisms. This is an
easy, and fairly well-known, consequence of the Bing Shrinking Criterion. (On
the other hand, I do not see any easy proof which does not use the BSC.)

\indent We proceed with the proof. If the nondegeneracy set of the given map
$f : M \rightarrow Q$ were compact, we would argue as follows. Given $\epsilon> 0$, let $\{
U_{i}|1 \leq i \leq p\}$ be a finite cover of the 0-dimension image of
nondeg$(f)$ by disjoint open subsets of $Q$ of diameter $< \epsilon$. Fixing
$i$, let $\{ N_{j}|1 \leq j <\infty\}$ be a strictly decreasing sequence of
compact, not-necessarily-connected manifold neighborhoods of nondeg$(f)
\cap f^{-1}(U_{i})$ in $f^{-1}(U_{i})$, such that $\cap^{\infty}_{j=1} N_{j} =
\mbox{nondeg}(f) \cap f^{-1}(U_{i})$, and such that each component of each
$N_{j}$ is null-homotopic in $N_{j-1}$. (Read $f^{-1}(U_{i})$ for $N_{0}$
here.) (The $N_{j^{\prime}s}$ can be chosen to be manifolds because without
loss $f^{-1}(U_{i})$ is a PL manifold, since each point inverse of $f$, being
cellular by say \cite{Mc}, has a PL manifold neighborhood.)

To connect the $N_{j}$'s, we start with $j = 1$, and proceed in increasing
order of the $j$'s, joining the components of $N_{j}$ together by tubes in
$int\; N^{\#}_{j-1}$ (let $N^{\#}_{0}$ be $P^{-1}(U_{i})$ here), to get a
compact connected manifold $N^{\#}_{j}$ which is null-homotopic in
$N^{\#}_{j-1}$. Then let $Y_{i} = \cap^{\infty}_{j=1}N^{\#}_{j}$, which is a
cell-like set containing nondeg$(f) \cap f^{-1}(U_{i})$. We can assume $Y_{i}$
has codemension $\geq3$, because the connecting operation can be done
carefully so that, for example, $Y_{i} - \mbox{nondeg}(f)$ has countably many
components, each a locally flatly embedded interval. Let $Q_{\#}$ be the
quotient space $M/\{Y_{i}|1 \leq i \leq p\}$. Then the quotient map $g : Q
\rightarrow Q_{\#}$ serves as the map $g$ in the preceding paragraph.

In the general case, when the nondegeneracy set of the given map $f : M
\rightarrow Q$ is not compact, but only $\sigma$-compact, one essentially does
a countable number of connecting operations as above, first for the
1-nondegeneracy set of $f$, then for the 1/2-nondegeneracy set of $f$, etc.,
where the $\epsilon$-nongeneracy set of $f$ is the compact set
$\cup\{f^{-1}(q)|\mbox{diam} f^{-1}(q) \geq\epsilon, q \epsilon Q\}$. But some
care, and explanation, is required. First, a little care is necessary to
ensure that the new intervals, which are introduced to connect nontrivial
point-inverses of $f$, miss all of the original nontrivial point-inverses.
This is easily done, using the codemension hypothesis. The second point,
requiring explanation, is more fundamental. One can do the connecting
operation on the 1-nondegeneracy set, then on the 1/2-nondegeneracy set
\underline{minus} the new 1-nondegeneracy set, then on the 1/3-nondegeneracy
set \underline{minus} the new 1/2-nondegeneracy set, etc., and this will
produce a countable upper semi-continuous decomposition, but it may not be
null, since at the second stage one may actually be producing a countable
number of new nondegenerate point-inverses of diameter $\geq1/2$. One way to
get around this is, when working on the 1/2-nondegeneracy set, to allow the
finitely many already-constructed new 1-nondegenerate point-inverses to
enlarge a little, by connecting to them the 1/2-nondegenerate point-inverses
which are sufficiently close. This way one can arrange to have only finitely
many 1/2-nondegenerate point-inverses at the end of the second stage. At the
next stage, one again allows an arbitrarily small enlargement of the already
constructed 1/2-nondegenerate point-inverses, so as to wind up with only a
finite number of 1/3-nondegenerate point-inverses. Since the amount of
enlarging is arbitrarily small at each stage, one can arrange in the limit
that the nondegenerate point-inverses be cell-like and codemension $\geq3$, as well as countable and null.

\indent The following rigorization of the above proof turns out to be a little
bit different in detail, but is the same in spirit, and has the virtue of
relative brevity. Basically, the idea is to intertwine the choice of the
connecting tubes with the choice of the manifold neighborhood sequence.

\indent Let $N_{1}$ be a not-necessarily-connected compact manifold
neighborhood of 1-nondeg$(f)$, so small that each component of $N_{1}$ $(i)$
lies in the 1-neighborhood of some point-inverse in 1-nondeg$(f)$, and $(ii)$
has image in $Q$ of diameter $< \epsilon$ (where $\epsilon> 0$ is given at the
start, as explained earlier). Let $N_{2} = N_{2,a} \cup N_{2,b}$ be a
not-necessarily-connected compact manifold neighborhood of 1/2-nondeg$(f)$,
where $N_{2,a}$ and $N_{2,b}$ are disjoint and each is a union of components
of $N_{2}$, so small that each component of $N_{2}$ $(i)$ lies in the
1/2-neighborhood of some point-inverse in 1/2-nondeg$(f)$, and $(ii)$ has image
in $Q$ of diameter $< \epsilon$, and so that $(iii)$ $N_{2,a}$ is a neighborhood
of 1-nondeg$(f)$, $(iv)$ each component of $N_{2,a}$ lies in, and is null
homotopic in, some component of $N_{1}$ (but the components of $N_{2,b}$ may
have no relation to $N_{1}$ at all), and $(v)$ each component of $N_{2,b}$ has
diameter $< 1$. One way to do all of this choosing is first to find a
partitioning $f$(1/2-nondeg$(f)) = C_{2,a} \cup C_{2}$ of the image of the
1/2-nondegeneracy set of $f$ into two disjoint closed (0-dimensional) subsets
$C_{2,a}$ and $C_{2,b}$ such that 1-nondeg$(f) \subset f^{-1}(C_{2,a})
\subset\; int\; N_{1}$. Then choose in the quotient $Q$ two disjoint open
neighborhoods $U_{2,a}$ of $C_{2,a}$ and $U_{2,b}$ of $C_{2,b}$, so small that
the preimages $f^{-1}(U_{2,a})$ and $f^{-1}(U_{2,b})$ and their components
satisfy conditions $(i) - (v)$ above. Now let $N_{2,a}$ be any compact manifold
neighborhood of 1/2-nondeg$(f) \cap f^{-1}(U_{2,a})$ in $f^{-1}(U_{2,a})$, and
likewise choose $N_{2,b}$.

\indent At this time, we tube together the various components of $N_{2,a}$
which lie in a common component of $N_{1}$ (but the components of $N_{2,b}$)
are not tubed together). We would like these tubes to miss $N_{2,b}$; a priori
that may not be possible, because $N_{2,b}$ may disconnect some of the
components of $N_{1}$. So what we do is first to choose the various connecting
tubes in $N_{1} - \; int\; N_{2,a}$, so thin and well-positioned that they
miss 1/2-nondeg$(f)$, but possibly may intersect $N_{2,b}$, and then we throw
away from $N_{2,b}$ a small neighborhood of the intersection of the tubes with
$N_{2,b}$, producing a smaller manifold $N^{\prime}_{2,b}$ which can take the
place of $N_{2,b}$. Finally, let $N^{*}_{2}$ denote the union of $N^{\prime
}_{2,b}$ and the tubed-together components of $N_{2,a}$. So $N^{*}_{2}$ has
one component for each component of $N_{1}$, and one component for each
component of $N^{\prime}_{2,b}$.

\indent This process is now repeated, to construct $N^{*}_{3}, N^{*}_{4}$,
etc. To save words, let it suffice to say that, in order to construct
$N^{*}_{i+1}$ given $N^{*}_{i}$, the above procedure works word-for-word,
after making the substitution $N^{*}_{i}$ for $N_{1}$; $N_{i+1}$, $N_{i+1,a}$
and $N_{i+1,b}$ for $N_{2}, N_{2,a}$ and $N_{2,b}$ and likewise for the $C$'s
and $U$'s; $N^{*}_{i+1}$ for $N^{*}_{2}$; and elsewhere $1/i$ for $1$ and
$1/(i+1)$ for $1/2$. (The very first substitution listed is the single exception
to the general theme $1 \rightarrow i$ and $2 \rightarrow i + 1$.)

\indent As part of the construction, one is obtaining at each stage an
\underline{injective} correspondence $\alpha_{i} :$ components of $N^{*}_{i}
\rightarrow$ components of $N^{*}_{i+1}$, such that for any component $P$ of
$N^{*}_{i}$, $\alpha_{i}(P)$ is null-homotopic in $P$. Hence for each $P,
\cap^{\infty}_{j=i}\alpha_{j}(\alpha_{j-1}(\ldots(\alpha_{i}(P))\ldots))$ is a
cell-like set, $P^{*}$ say, of codemension $\geq3$ by the usual additional
case. (This latter claim uses $(i)$, to conclude that the set $P^{*}$ minus the
1-demensional connecting intervals in $P^{*}$ lies in nondeg$(f)$, hence the
codemension $\geq3$). Then the null collection $\{P^{*}\}$ consisting of all
of these cell-like sets, exactly one for each component of each difference
manifold $N^{*}_{i} - \cup\alpha_{i-1}$ (components of $N^{*}_{i-1}$), $i
\geq2$ (let $N^{*}_{1} = \emptyset$ here), is the desired null collection.
This completes the proof that Countable Shrinking Theorem $\Rightarrow$
0-Dimensional Shrinking Theorem.
\newline%
\newline\textbf{\textit{Proof of the Countable
Shrinking Theorem.}} Let $\{Y_{1}, Y_{2}, \ldots\}$ be the countable null
collection of disjoint codimension $\geq3$ cell-like sets in $M$, which are
the nontrivial point-inverses of $f$. Our goal is to prove: 
\newline%
\newline\textbf{Shrinking Lemma:} Given any $Y_{i}$ and any $\epsilon> 0$,
there is a neighborhood $U$ of $Y_{i}$, $U \subset N_{\epsilon}(Y_{i})$, and
there is a homeomorphism $h : M \rightarrow M$, supported in $U$, such that
for any $Y_{j}$, if $h(Y_{j}) \cap U \neq\emptyset$, then $diam$ $h(Y_{j}) <
\epsilon$.
\newline\newline \textbf{Technical Note.} The reason for writing
$h(Y_{j})\cap U \neq\emptyset$ instead of the equivalent $Y_{j} \cap U
\neq\emptyset$, is to make this statement more nearly resemble that of the
a-Shrinking Lemma, below.
\newline
\newline\indent Given this Lemma, it is an easy matter to show that the Bing Shrinking
Criterion holds for the given decomposition $f: M \rightarrow M/\{Y_{i}\} = Q$
as in \cite{Bi1}. For given $\epsilon> 0$, the BSC is that there
exist a homeomorphism $h : M \rightarrow M$ such that $(i)$ $dist (fh,f) <
\epsilon$, and $(ii)$ for each $Y_{i}$, $diam$ $h(Y_{i}) < \epsilon$. Such an $h$
can be gotten by applying the above Shrinking Lemma to sufficiently small
disjoint neighborhoods of the finitely many $Y_{i}$'s which initially have
diameter $\geq\epsilon$.

\indent As an introduction to the proof of the Shrinking Lemma, we consider
the trivial demension range case, where each $Y_{i}$ satisfies $2\,dem\,Y_{i} +
2 \leq m$. In this case, to prove the Lemma for $Y_{i}$ say, one starts by
embedding $cY_{i}$ (the cone on $Y_{i}$) in $N_{\epsilon}(Y_{i})$, extending
the given embedding $Y_{i} \hookrightarrow N_{\epsilon}(Y_{i})$ of its base
$Y_{i}$, so that $dem\,cY_{i} \leq dem\,Y_{i} + 1$. (If $2\,dem\,Y_{i} + 3 \leq
m$, this is a classical result of Menger-N\"{o}beling \cite{Hu-Wa}; if $2\,dem\,Y_{i} +
2 = m$, see the next paragraph for how to do something just as good.) Then, by
general position applied a countable number of times, one can ambient isotope
$cY_{i}$ an arbitrarily small amount, keeping its base $Y_{i}$ fixed, to make
$cY_{i}$ disjoint from all of the other $Y_{j}$'s. Thus, one now has a
guideway in $M$ minus the other $Y_{j}$'s, namely $cY_{i}$, along which to
shrink $Y_{i}$. Of course, the other $Y_{j}$'s may converge closer and closer
to $cY_{i}$, but as they do so, they must get smaller and smaller, by their
nullity. The exact method of constructing the shrinking homeomorphism $h$ at
this point is modeled on Bing's original work, and is discussed fully in the
following paragraphs.

\indent Given a cellular set $Y$ in $M$ (for example, a codemension $\geq3$
cell-like set) and a neighborhood $U$ of $Y$, we can find a coordinate
chart $\mathbb{R}^{m} \hookrightarrow M$ of $M$ lying in $U$ such that $Y \subset int
\; B^{m}$, where $B^{m}$ is the unit ball in the coordinate chart, and also
origin $0 \notin Y$. (We hope the reader will be able to tolerate
this frequently-occurring abuse of notation, namely out not labeling the
embedding $R^{m} \hookrightarrow M$. This will leave upcoming expressions much
less cluttered.) Let $Y_{1} \subset\partial B^{m}$ denote the image of $Y$,
projected from the origin. (Similarly, for upcoming use, let $Y_{r}
\subset\partial\,rB^{m}$ denote the projected image of $Y$, where $rB^{m}$ is
the ball of radius $r > 0$.)
\\
\newline\textbf{Lemma.} \textit{By an
arbitrarily small perturbation of the coordinate chart embedding, we can
arrange that dem $Y_{1} \leq$ dem $Y$.}
\\
\\
\indent To minimize ambiguity here, we emphasize that $Y$, being regarded as a
subset of $M$, is left fixed; it is only the coordinate chart embedding $\mathbb{R}^{m}
\hookrightarrow M$ that is being changed, in order to change what (the
preimage of) $Y$ looks like in $\mathbb{R}^{m}$. Note that some perturbing may be
necessary; even a tame cantor set $Y$ in int $B^{m} - 0$ can have its
projection-image $Y_{1}$ be all of $\partial B^{m}$.
\\
 \newline\textbf{\textit{Proof of Lemma:}} If $Y$ were a locally flatly (on each open stratum)
embedded polyhedron in int $B^{m}$, one could achieve the Lemma by making $Y$
PL embedded in $B^{m}$ (PL on each open stratum would suffice). The general
proof is the demension theory analogue of this. Actually, it is quickest to
think in terms of complements. Let $K = L^{m-y-2} \cup\infty$ be a countable
union of $(m-y-2)$-sphere on $\partial B^{m} \approx\mathbb{R}^{m-1} \cup\infty$,
where $y = \dim Y$, gotten from the countable union of $(m-y-2)$-planes in
$\mathbb{R}^{m-1}$ consisting of all points of $\mathbb{R}^{m-1}$ having at least $y + 1$
coordinates rational. (Thus $\partial B^{m} - K \approx\mathbb{R}^{m-1} - L$ is
N\"{o}beling's $y$-dimensional space.) Any compact subset of $\partial
B^{m}-K$ has demension $\leq y$. To arrange that $Y_{1} \subset\partial B^{m}
- K$, one applies general position, perturbing the coordinate chart
embedding on int$B^{m}$ to make $cK$ (which has dimension $m-y-1$) disjoint
from $Y$. $\blacksquare$

\indent Assuming then that dem $Y_{1} \leq$ dem $Y$, let $X = cY_{1}$, that
is, $Y_{1}$ coned to the origin in $B^{m}$. By construction, dem $X \leq$ dem
$Y +1$. (For upcoming use let $rX = cY_{r}$; it has the same demension as $X$).

\indent To such an $X$, we associate a \underline{fixed sequence of special
compact neighborhoods} $N_{1}\supset int\;N_{1}\supset
N_{2}\supset\ldots$ of $X$ in $M$, with $X=\cap_{\ell=1}^{\infty}N_{\ell}$,
constructed as follows. Let $O_{1}\supset int\;O_{1}\supset O_{2}\supset
\ldots$ be any strictly decreasing sequence of compact neighborhoods of
$Y_{1}$ in $\partial B^{m}$ such that $Y_{1}=\cap_{\ell=1}^{\infty}O_{\ell}%
$. Define $N_{\ell}=c(1+1/\ell)O_{\ell}\cup1/\ell B^{m}$, that is, $N_{\ell
}$ is the cone (to the origin) on the projection image $(1+1/\ell)O_{\ell
}\subset\partial(1+1/\ell)B^{m}$, together with the $m$-ball of radius $1\ell$. (See Figure 1).

\begin{figure}[ht]
\centerline{
\includegraphics{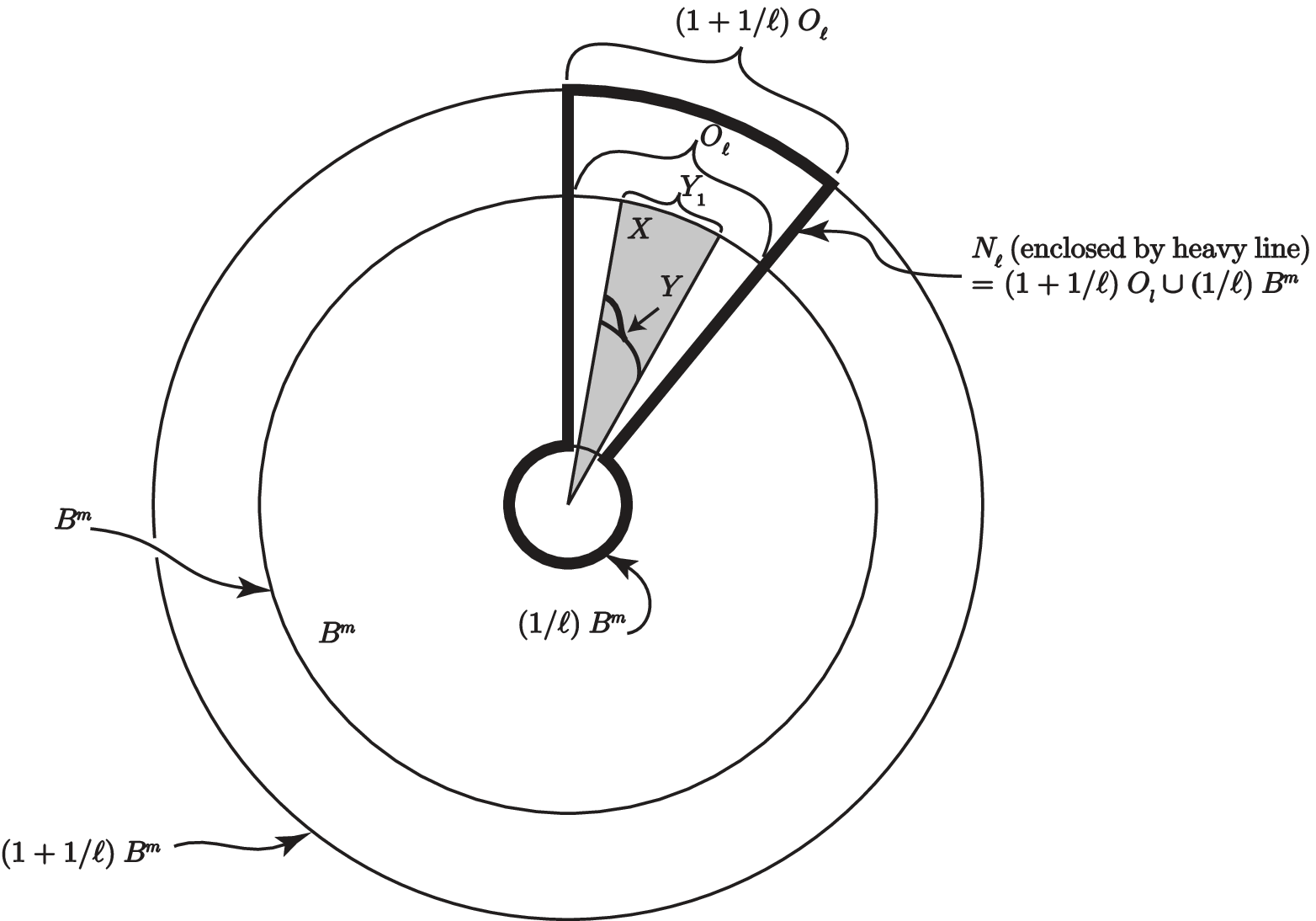}
}
\caption{}
\end{figure}

We emphasize that given $Y$, once its associated $X$ and its neighborhood
sequence $\{ N_{\ell}\}$ have been constructed, they are fixed entities for the duration of
the proof (except that perhaps certain \underline{arbitrarily small} general
positioning ambient isotopies may be applied to them).

\renewcommand{\thefootnote}
{\ensuremath{\fnsymbol{footnote}}}
\setcounter{footnote}{1}

\indent We now describe the fundamental squeezing homeomorphism that will be
used throughout the proof. Given $X$ and $\{ N_{\ell}\}$ as above, and given a
finite indexing subsequence $\lambda: \{1, \ldots, p\} \rightarrow\{p, p+1,
\ldots\}$ of length $p$, where $p \geq2$ is arbitrary, say $\lambda: p
\leq\lambda(1) < \lambda(2) < \ldots< \lambda(p)$, we define a specific
homeomorphism $h_{\lambda}: M \rightarrow M$ which will have the following
properties:\\[6pt]
$(i)$ $h_{\lambda}$ is fixed off of $N_{\lambda(1)}$
and on $1 / \lambda(1)B^{m}$,\\[6pt]
$(ii)$ $h_{\lambda}$ is radial and
inward-moving in the coordinate chart structure, i.e., each radial line
segment in $2B^{m}$ is carried onto itself by $h_{\lambda}$, and each point,
if moved, is moved toward the origin,\\[6pt]
$(iii)$ $h_{\lambda}(X) \subset2/p\,B^{m}$, and\\[6pt]
$(iv)$ for any connected subset $C$ of $M$ which intersects at most a
single $fr\,N_{\lambda(j)}$, the radial-height of $C$ (defined below) is
increased by at most $4/p$\footnote{Actually, $2/p$ will work here, but my crude analysis does not yield that.} under $h_{\lambda}$. Hence, if such a subset $C$ lying in the coordinate patch $\mathbb{R}^{m} \hookrightarrow M$ has euclidean-metric diameter $\leq\eta$, then $h_{\lambda}(C)$ has euclidean-metric diameter $\leq\eta+ 4/p$.
\\\\
\indent The \underline{radial height} of $C$ is the length of the projection
of $C \cap\mathbb{R}^{m}$ onto the radius coordinate $[0, \infty)$ in the coordinate
chart $\mathbb{R}^{m} \hookrightarrow M$. (Since $C$ is connected, this
projection-image of $C \cap \mathbb{R}^{m}$ is an interval. To be technically complete,
let us agree that in the special case when this image is the interval $[c,
\infty)$, then $(iv)$ above means that the projection-image of $h(C) \cap \mathbb{R}^{m}$
lies in $[c - 4/p, \infty)$.)

\indent On each radial line segment of $2B^{m}$, $h_{\lambda}$ will be
piecewise linear, with at most $p+1$ \textquotedblleft
breaks\textquotedblright\ (changes in derivative); they will be at the
(source) levels $\partial(1+1/ \lambda(j))B^{m}$, $1\leq j\leq p$, and at
$\partial(1/\lambda(1))B^{m}$. Let these $p+1$ decreasing radii be
$s(1),\ldots,s(p+1)$. That is, for each $j,1\leq j\leq p$, let
$s(j)=1+1/\lambda(j)$, and let $s(p+1)=1/\lambda(1)$. Let $t(1),\ldots,t(p+1)$
be the $p+1$ equally spaced numbers decreasing from $1+1/\ell(1)$ to
$1/\ell(1)$; inclusive. That is, if each $j,1\leq j\leq p+1$, let
$t(j)=1+1/\lambda(1)-(j-1)/p$. The general idea is to have $h_{\lambda}$ carry
the $s$-levels to the $t$-levels. We start by defining $h_{\lambda}$ on the
\textquotedblleft outer end\textquotedblright\ of each $N_{\lambda(j)}$,
namely on $s(j)O_{\lambda(j)}$, so that $h_{\lambda}(s(j)O_{\lambda
(j)}=t(j)O_{\lambda(j)}$. (Necessarily then the image $h_{\lambda}%
(N_{\lambda(j)})$ must be $ct(j)O_{\lambda(j)}\cup1/\lambda(1)B^{m}$.) Next,
define $h_{\lambda}$ on the cylindrical sleeve $[s(p+1),s(1)](O_{\lambda
(j)}-int\;O_{\lambda(j+1)}), 1\leq j\leq p$, (here $[a,b]O$ denotes
$cl(bO-aO)$, and $O_{\lambda(p+1)}=\phi)$ one at a time, in order of
increasing $j$. At the $j$th stage, $h_{\lambda}$ has already been defined on
$[s(p+1),s(1)]fr\;O_{\lambda(j)}$, and has there $j+1$ breaks (or no breaks,
if $j=1$), at the $j+1$ (source) levels $s(1)fr\;O_{\lambda(j)},\ldots
,s(j)fr\;O_{\lambda(j)},s(p+1)fr\;O_{\lambda(j)}$, which are taken
respectively by $h_{\lambda}$ to the $j+1$ (target) levels $t(1)fr\;O_{\lambda
(j)}\ldots,t(j)fr\;O_{\lambda(j)},t(p+1)fr\;O_{\lambda(j)}$. Define
\begin{align*}
h_{\lambda}|[s(j),&s(1)](O_{\lambda(j)} - int\;O_{\lambda(j+1)})  =\\[6pt]
& h_{\lambda}|[s(j),s(1)]fr\;O_{\lambda(j)}\times\;\mbox{identity}(O_{\lambda
(j)}-int\;O_{\lambda(j+1)}),
\end{align*}
where the meaning of this expression should be clear. In order to define
$h_{\lambda}$ on the remaining region $[s(p+1),s(j)](O_{\lambda(j)}%
-int\;O_{\lambda(j+1)})$ to complete this stage of the definition, we must
choose a Urysohn function to tell us where to send the level
$s(j+1)(O_{\lambda(j)}-int\;O_{\lambda(j+1)})$. Let
\begin{align*}
\phi:s(j+1)(O_{\lambda(j)}-int\;O_{\lambda(j+1)})\longrightarrow\lbrack
t(j+1),t(j)]
\end{align*}
be a map such that $\phi|s(j+1)fr\;O_{\lambda(j)}=$ the radius value of $h_{\lambda
}(s(j+1)fr\;O_{\lambda(j)})$ (which one may compute) and $\phi
|s(j+1)fr\;O_{\lambda(j+1)}=t(j+1)$. Then define $h_{\lambda}$ on the level
$s(j+1)(O_{\lambda(j)}-int\;O_{\lambda(j+1)})$ by $h_{\lambda}(s(j+1)x)=\phi
(x)\cdot x$ for each $x\in O_{\lambda(j)}-int\;O_{\lambda(j+1)}$. Finally
extend $h_{\lambda}$ in linear fashion over each radial line segment on the
regions $[s(p+1),s(j+1)](O_{\lambda(j)}-int\;O_{\lambda(j+1)})$ and
$[s(j+1),s(j)](O_{\lambda(j)}-int\;O_{\lambda(j+1)})$. This completes the
definition of $h_{\lambda}$.

\begin{figure}
\begin{tabular}{c}
\includegraphics{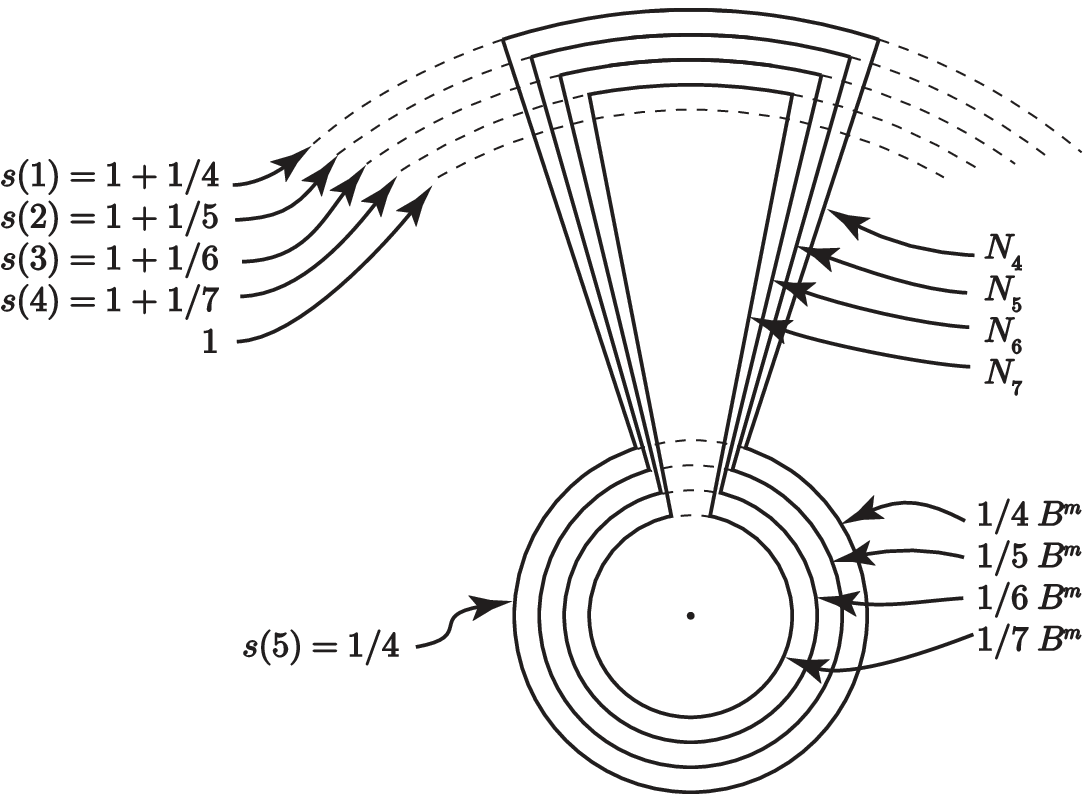}\\
Before applying $h_\lambda$\\
\bigskip{}\\
\includegraphics{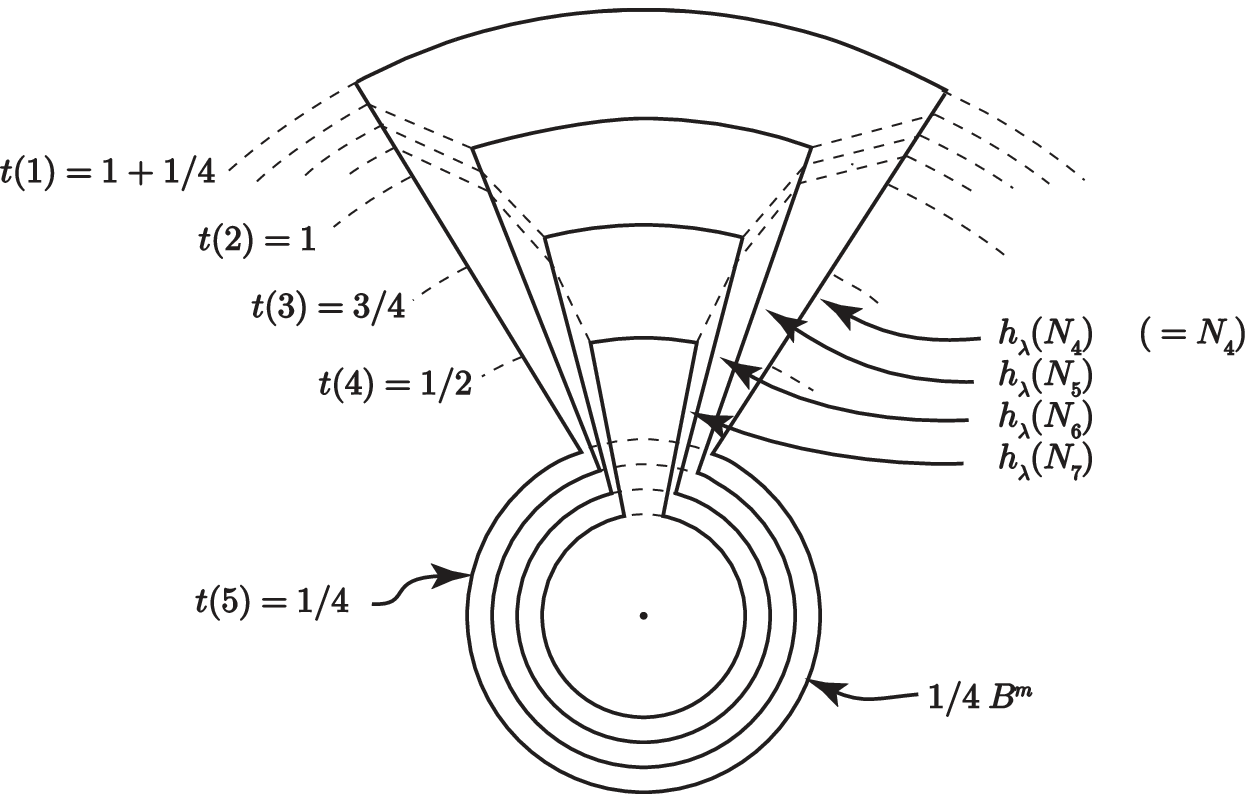}\\
After applying $h_\lambda$
\end{tabular}
\caption{The squeezing homeomorphism $h_\lambda$.  Pictured here is the case $p=4=\lambda(1)<\dots<\lambda(4)=7$.}
\end{figure}

\begin{figure}[ht]
\centerline{
\includegraphics{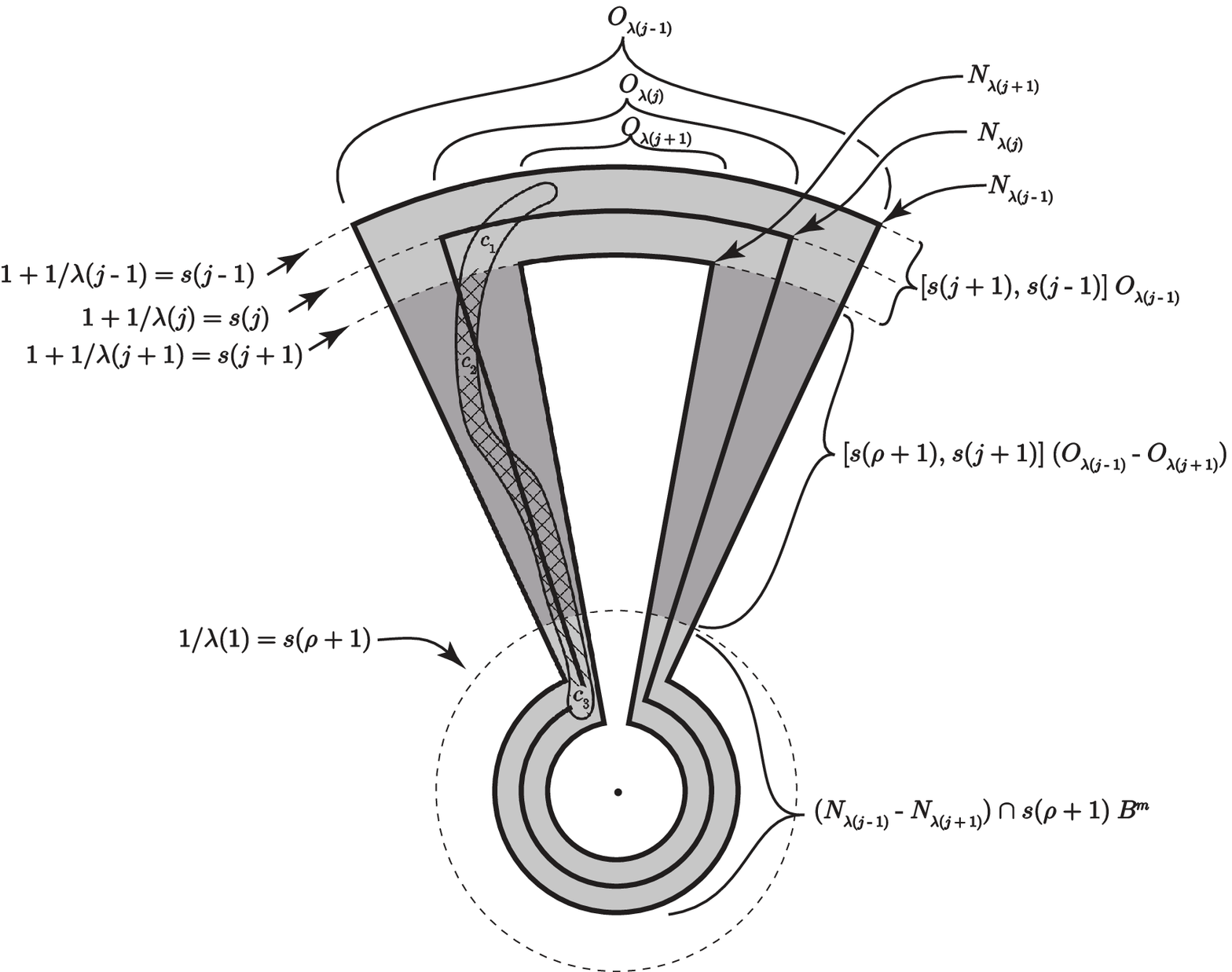}
}
\caption{}
\end{figure}

\indent The only nontrivial property of $h_{\lambda}$ to verify is $(iv)$, and
that can be understood by looking at Figures 2 and 3. The $h_{\lambda}$-image of an arbitrary connected set $C\subset N_{\lambda (j-1)}-N_{\lambda(j+1)}$ can be analyzed by breaking $C$ into three pieces. First there is $C_{1}\equiv C\cap\lbrack s(j+1),s(j-1)]O_{\lambda(j-1)}$,
whose image under $h_{\lambda}$ has radial height $\leq2/p$. Then there is
$C_{2}\equiv C\cap\lbrack s(p+1),s(j+1)](O_{\lambda(j-1)}-O_{\lambda(j+1)})$.
On each radial interval of $[s(p+1),s(j+1)](O_{\lambda(j-1)}-O_{\lambda
(j+1)})$ $h_{\lambda}$ is \underline{linear} (and compressing), fixed on the
inner end of the interval, while all of the outer ends (i.e., $s(j+1)(O_{\lambda(j-1)}-O_{\lambda(j+1)})$ have their images under $h_{\lambda}$ in the region $[t(j+1),t(j-1)]B^{m}$. So the radial height of
$C_{2}$ is increased by at most $2/p$, which is the difference between $t(j-1)$ and $t(j+1)$. Finally, there is $C_{3}\equiv C\cap s(p+1)B^{m}$, which is left fixed by $h_{\lambda}$. Combining these facts, one obtains property $(iv)$. This completes our discussion of the standard squeezing
homeomorphism $h_{\lambda}$.

\indent We return for a moment to the trivial demension range case, when $2\;
dem\;Y_{i} + 2 \leq m$ for each $i$, to illustrate exactly how the
above-constructed shrinking homeomorphism will be called into play. Having
fixed $i$ and $\epsilon> 0$, earlier we said to nicely embed $cY_{i}$ in
$N_{\epsilon}(Y_{i})$; now instead we find a cone $X_{i} = cY_{i,1}$ containing
$Y_{i}$, constructed as above so that $dem\: X_{i} \leq dem\: Y_{i} + 1$ by
working in a coordinate chart of $M$ which lies in $N_{\epsilon}(Y_{i})$. As
earlier, $X_{i}$ can be general positioned to intersect none of the other
$Y_{j} ^{\prime}s$. Let $\{N_{\ell}\}$ be a fixed sequence of neighborhoods of
$X_{i}$, constructed as above. The goal now is to choose $p \geq2$ and a
subsequence $\lambda: p \leq\lambda(1) < \lambda(2) < \ldots< \lambda(p)$ so
that the associated squeezing homeomorphism $h_{\lambda}$ satisfies the
conclusion of the Shrinking Lemma. The precise way to choose this subsequence
$\lambda$ is explained in the proof of the following Proposition. We emphasize
that in this Proposition, all distances are measured in the given metric on
the manifold $M$.
\\
\newline\noindent\textbf{Squeezing Proposition.} Suppose $X
\subset M$ is any compact cone lying in a coordinate chart of $M$ as described
above (that is, $X = c(X \cap\partial B^{m})$, where $B^{m}$ is the standard
$m$-ball in the coordinate chart), and suppose $\{N_{\ell}|1 \leq\ell<
\infty\}$ is a sequence of compact neighborhoods of $X$ as described above.
Then given $\epsilon> 0$ there exists $\delta> 0$ and an integer $p \geq2$
such that for any subsequence $\lambda: p \leq\lambda(1) < \lambda(2) <
\ldots< \lambda(p)$ of integers, the squeezing homeomorphism $h_{\lambda}$
(described above) has the following properties:\\[6pt]
\indent (1) $diam\: h_{\lambda}(X) < \epsilon$, and\\[6pt]
\indent (2) for any connected subset $C$ of $M$ such that $diam\: C < \delta$
and $C$ intersects at most one of the sets $fr\: N_{\lambda(j)}$, $1 \leq j \leq
p$, one has that $diam\: h_{\lambda}(C) < \epsilon.$
\\\\
\textbf{\textit{Proof.}} The proof of properties (1) and (2) rests on the fact that the
two metrics on $\mathbb{R}^{m} \hookrightarrow M$, the one induced from the $M$-metric
and the other being the standard euclidean metric, are equivalent. Since
$h_{\lambda}$ has support in $2B^{m}$, let us assume without loss that the set
$C$ of (2) lies in $3B^{m}$. Given $\epsilon> 0$, let $\epsilon^{\prime}> 0$
be such that any subset of $3 B^{m}$ having euclidean-metric-diameter $<
\epsilon^{\prime}$ has $M$-metric-diameter $< \epsilon$. Let $\delta^{\prime}>
0$ and $p \geq2$ be such that $\delta^{\prime}+ 4/p < \epsilon^{\prime}$. Finally, let $\delta> 0$ be such that any subset of $3B^{m}$ having $M$-metric-diameter $< \delta$ has euclidean-metric-diameter $< \delta^{\prime}$. Now given any sequence $\lambda: p \leq\lambda(1) < \ldots<
\lambda(p)$, and given any connected set $C$ as in the proposition, $C \subset3B^{m}$, then the euclidean-metric-diameter of $C$ is $< \delta^{\prime}$, and hence the euclidean-metric-diameter of $h_{\lambda}(C)$ is $\delta^{\prime}+ 4/p < \epsilon^{\prime}$, and hence the $M$-metric-diameter
of $h_{\lambda}(C)$ is $< \epsilon$. $\blacksquare$
\\\\
\indent In order to prove the general codemension 3 case of this theorem, we
do an iterated general positioning operation, just as in the proof of
codemension $\geq3$ engulfing. The inductive hypothesis is provided by the
following Lemma; the Shrinking Lemma above can be thought of as the $A=M$ case
of this Lemma. 
\\\\
\textbf{$\boldsymbol{a}$-Shrinking Lemma $(-1 \leq a \leq
m-2)$.} \textit{Suppose $A$ is a closed subset of $M$, with dem $A \leq a$.
Given any $Y_{i}$ and any $\epsilon> 0$, there is a neighborhood $U$ of
$Y_{i}, U \subset N_{\epsilon}(Y_{i})$, and there is a homeomorphism $h : M
\rightarrow M$, supported in $U$, such that for any $Y_{j}$, if $h(Y_{j}) \cap
U \cap A \neq\varnothing$, then diam $h(Y_{j}) < \epsilon$.}
\\\\
\indent This will be proved by induction on increasing $a$. But first we illustrate
the general idea by establishing that $(m-2)$-Shrinking Lemma $\Rightarrow$
Shrinking Lemma. We point out, for the reader who would like to gain familiarity with the entire proof a step at a time (as I did), that\\[6pt]
\indent (1) In the trivial dimension range case, when $2y + 2 \leq m$, where
$y = \max\{dem\: Y_{i}\}$, one is in effect only using the a-Shrinking Lemma for
$a \leq y + 1$, whose proof is a trivial general position argument, already
used, together with the proof that $(y+1)$-Shrinking Lemma $\Rightarrow$
Shrinking Lemma, as below.\\[6pt]
\indent (2) in the ``metastable'' dimension range case, when $3y + 4 \leq2m$
($y$ as above), one is in effect only using the a-Shrinking Lemma for $a \leq
y + 1$, whose proof in turn only uses the a-Shrinking Lemma for $a \leq2y + 3
- m$, whose proof is the aforementioned trivial general position argument,
together with the proof that $(y + 1)$-Shrinking Lemma $\Rightarrow$ Shrinking
Lemma, as below.
\\\\
\textbf{\textit{Proof that $\boldsymbol{(m-2)}$-Shrinking Lemma $\Rightarrow$ Shrinking Lemma.}} 
Denote by $Y_{0}$ the given $Y_{i}$ in the Shrinking Lemma. Given $\epsilon> 0$, 
let $U_{0} \subset N_{\epsilon}(Y_{0})$ be a saturated neighborhood of $Y_{0}$ (\underline{saturated} meaning that if any
$Y_{j}$ intersects $U_{0}$, it lies in $U_{0}$), so small that if $Y_{j} \subset
U_{0}$, $Y_{j} \neq Y_{0}$, then $diam \; Y_{j} < \epsilon$. In $U_{0}$,
choose a coordinate patch of $M$ containing $Y_{0}$, and construct there in
the manner explained earlier a cone $X_{0}$ containing $Y_{0}$, with dem
$X_{0} \leq dem \; Y_{0} + 1 \leq m - 2$, and a special neighborhood basis $\{
N_{\ell}\}$ of $X_{0}$. Let $\delta> 0$ and $p \geq2$ be as provided by the
Proposition , for this data (without loss $\delta< \epsilon)$. Our goal is to
move off of $X_{0}$ those intersecting $Y_{j}$'s (other than $Y_{0}$) which
are too big, by using the $(m-2)$-Shrinking Proposition, leaving behind to
intersect $X_{0}$ only $Y_{j}$-images of size $< \delta$. Then we iterate this
operation $p-1$ more times in order to achieve the desired insulation of
$X_{0}$ from $M - N_{\lambda(1)}$. We remark now that, even though these
various moves in the successive stages may have overlapping supports, their
composition will not stretch any $Y_{j}$ to have diameter $\geq\epsilon$.

\indent To start, let $\mathscr{Y}_{1}$ be the finite subcollection of members of $\mathscr{Y} -
\{ Y_{0} \}$ which intersect $X_{0}$ and have diameter $\geq\delta$. (Here $\mathscr{Y}$
denotes the entire collection $\{Y_{i}\}$.) Choose a collection $\mathscr{U}_{1}$ of
disjoint open saturated neighborhoods of the members of $\mathscr{Y}_{1}$, each member
of $\mathscr{U}_{1}$ having diameter $< \epsilon$ and lying in $U_{0} - Y_{0}$. For each
member $U$ of $\mathscr{U}_{1}$, apply the $(m-2)$-Shrinking Lemma, with $\epsilon
$-value $\min\{\delta, dist\,(\cup\, \mathscr{Y}_{1}, M - \cup\, \mathscr{U}_{1})\}$, to find a
homeomorphism $h_{U}$, supported in $U$, such that the $h_{U}$-image of any
member of $\mathscr{Y}$ lying in $U$ and intersecting $X_{0}$ has diameter $< \delta$.
Letting $H_{1}$ be the composition of these $h_{U}$'s, $U \in \mathscr{U}_{1}$, it
follows that each member of $H_{1}(\mathscr{Y} - \{Y_{0}\})$ which intersects $X_{0}$
has diameter $< \delta$. So we can choose $\lambda(1) \geq p$ so large that
$N_{\lambda(1)} \subset U_{0}$ and each member of $H_{1}(\mathscr{Y} - \{Y_{0}\})$ which
intersects $N_{\lambda(1)}$ has diameter $< \delta$,

\indent From now on, the repeating steps are qualitatively the same, but they
are a little bit different from the just-completed first step. Let $\mathscr{Y}_{2}$ be
the finite subcollection of members of $H_{1}(\mathscr{Y} - \{Y_{0}\})$ which intersect
both $fr\,N_{\lambda(1)}$ and $X_{0}$. (Possibly $\mathscr{Y}_{2} \cap H_{1}(\mathscr{Y}_{1})
\neq\varnothing$; that is allowable.) Choose a collection $\mathscr{U}_{2}$ of disjoint open
$H_{1}(\mathscr{Y})$-saturated neighborhoods of the members of $\mathscr{Y}_{2}$, each having
diameter $< \delta$ (which is $< \epsilon$) and lying in $U_{0} - Y_{0}$. For
each member $U$ of $\mathscr{U}_{2}$, apply the $(m-2)$-Shrinking Lemma, with $\epsilon
$-value $\min\{dist\,(X_{0}, M-N_{\lambda(1)}),dist\,(\cup\,\mathscr{Y}_{2}, M - \cup\,
\mathscr{U}_{2})\}$, to find a homeomorphism $h_{U}$, supported in $U$, such that the
$h_{U}$-image of any member of $H_{1}(\mathscr{Y})$ which intersects $X_{0}$ necessarily
misses $fr\,N_{\lambda(1)}$. Letting $H_{2}$ be the composition of these
$h_{U}$'s, $U \in \mathscr{U}_{2}$, it follows that each member of $H_{2} H_{1}(\mathscr{Y} -
\{Y_{0} \})$ which intersects $N_{\lambda(1)}$ has diameter $< \delta$, and
also each member of $H_{2} H_{1}(\mathscr{Y})$ intersects at most one of $fr\,
N_{\lambda(1)}$ and $X_{0}$. We can now choose $\lambda(2) > \lambda(1)$ so
large that each member of $H_{2} H_{1}(\mathscr{Y})$ intersects at most one of $fr\,
N_{\lambda(1)}$ and $fr\,N_{\lambda(2)}$.

\indent In general, the argument goes as follows. (The following $k=2$ case was
done above). Given $k,\ 2\leq k\leq p$, suppose we have constructed a
homeomorphism $G_{h-1}(= H_{k-1}\circ\ldots\circ H_{1}) : M \rightarrow M$,
supported in $U_{0} - Y_{0}$, and a sequence $p \leq\lambda(1) < \ldots<
\lambda(k-1)$, with the properties:\\[6pt]
\indent  $(1_{k-1})$ each member of $G_{k-1}(\mathscr{Y} - \{ Y_{0} \})$ lying in $U_{0}$ has diameter $< \epsilon$, and each member which intersects $N_{\lambda(1)}$ has diameter $< \delta$, and\\[6pt]
\indent$(2_{k-1})$ each member of $G_{k-1}(\mathscr{Y})$ intersects at most one
of $fr\,N_{\lambda(1)}, \ldots, fr\,N_{\lambda(k-1)}$.\\[6pt]
We show how to construct the analogous $G_{k}$ and $\lambda(k)$. Let $\mathscr{Y}_{k}$ be the finite subcollection of members of $G_{k-1}(\mathscr{Y} - \{Y_{0}\})$ which intersect both $fr\,N_{\lambda(k-1)}$ and $X_{0}$. Choose a collection $\mathscr{U}_{k}$ of disjoint open
$G_{k-1}(\mathscr{Y})$-saturated neighborhoods of the members of $\mathscr{Y}_{k}$, each having
diameter $< \delta\:(< \epsilon)$ and lying in $N_{\lambda(k-2)}-Y_{0}$. (Let
$N_{\lambda(0)}$ be $U_{0}$ here.) For each member $U$ of $\mathscr{U}_{k}$, apply the
$(m-2)$-Shrinking Lemma, with $\epsilon$-value $\min\{dist(X_{0},
M-N_{\lambda(k-1)})$, $dist(\cup \mathscr{Y}_{k}, M - \cup \mathscr{U}_{k})\}$, to find a
homeomorphism $h_{U}$, supported in $U$, such that the $h_{U}$-image of
any member of $G_{k-1}(\mathscr{Y})$ which intersects $X_{0}$ necessarily misses $fr\,N_{\lambda(k-1)}$. Letting $H_{k}$ be the composition of these $h_{U}$'s, $U
\in \mathscr{U}_{k}$, and letting $G_{k} = H_{k} \circ G_{k-1}$, it follows that
$G_{k}$ satisfies properties $(1_{k})$ and $(2^{\prime}_{k})$, where
$(2^{\prime}_{k})$ is property $(2_{k})$ with $X_{0}$ in place of $fr\,
N_{\lambda(k)}$. To achieve $(2_{k})$, simply choose $\lambda(k) >
\lambda(k-1)$ so large that each member of $G_{k}(\mathscr{Y})$ intersects at most one
of $fr\,N_{\lambda(k-1)}$ and $fr\,N_{\lambda(k)}$.

\indent After constructing $G_{p}$ in this manner with properties $(1_{p})$
and $(2_{p})$, the final homeomorphism $h$ of the Shrinking Lemma is $h =
h_{\lambda}G_{p}$, where $h_{\lambda}$ is the squeezing homeomorphism
constructed earlier, corresponding to the finite sequence $\lambda: \lambda(1)
< \ldots< \lambda(p)$. It follows from the Squeezing Proposition that $h$,
which is supported in $U = U_{0}$, has the desired properties.

\indent This completes the proof that $(m-2)$-Shrinking Lemma $\Rightarrow$ Shrinking
Lemma.
\\\\
\textbf{\textit{Proof that $\boldsymbol{(a-1)}$-Shrinking Lemma $\Rightarrow$
$\boldsymbol{a}$-Shrinking Lemma, for $\boldsymbol{a \leq m -}$}} $\boldsymbol{2.}$ 
As the reader will recognize, this
proof is modelled on the preceding proof, but it is a wee bit more complicated.

\indent Let $A, Y_{i}$ and $\epsilon> 0$ be as in the hypothesis of the
$a$-Shrinking Lemma; as before let $Y_{0}$ denote this $Y_{i}$. Let $U_{0}
\subset N_{\epsilon}(Y_{0})$ be a saturated neighborhood of $Y_{0}$; so small
that if $Y_{j} \subset U_{0}$, $Y_{j} \neq Y_{0}$, then $diam \, Y_{j} <
\epsilon$. In $U_{0}$, choose a coordinate patch of $M, \mathbb{R}^{m} \hookrightarrow
M$, containing $Y_{0}$, and construct there in the manner explained earlier a
cone $X_{0}$ containing $Y_{0}$, with $dem\, X_{0} \leq dem\, Y_{0} + 1 \leq
m-2$. \underline{In addition}, we wish to construct $X_{0}$ to
be in general position with respect to $A$, in such a manner that $X_{0} \cap
A$ lies in a compact subcone $Z$ of $X_{0}$ (i.e. $Z = c(Z \cap c \,\partial B^{m})
\subset X_{0} = c(X_{0} \cap \partial B^{m})$) with the property that $dem \, Z \leq
dem \, A + dem \, X_{0} - m + 1$, hence $dem \, Z \leq a - 1$. One
way to do this is as follows. First construct $X_{0}$ in the manner described
earlier, without regard to $A$. Then, perturbing the coordinate chart
structure (i.e. perturbing the embedding $\mathbb{R}^{m} \hookrightarrow M$, as in the
earlier Lemma) an arbitrarily small amount (this time moving $X_{0}$ hence
$Y_{0}$, but regarding $A$ as being fixed), arrange that $A \cap B^{m} \subset
D = c D_{1}$, where $D$ is a subcone of $B^{m}$, $D_{1} = D \cap \partial B^{m}$, and
$dem \, D_{1} \leq dem \, A$ and hence $dem \, D \leq dem \, A + 1$. This is
done by the same argument used to construct $X_{0}$, by moving a certain
$\sigma$-compact cone $c K^{m-a-2}$ in $B^{m}$, of demension $m-a-1$, off of
$A$. Now, to construct $Z$, perturb the coordinate chart structure again
(again moving $X_{0}$ hence $Y_{0}$, but still thinking of $A$, and in
addition $D$, is fixed), by first isotoping $\partial B^{m}$ in itself to make the
set $X_{0} \cap\partial B^{m} (= X_{0,1}$ in earlier notation) in general
position with respect to $D_{1}$ in $\partial B^{m}$ (i.e., dem\,$X_{0} \cap\partial
B^{m} \cap D_{1} \leq dem\,X_{0} \cap \partial B^{m} + dem \, D_{1} - (m-1))$ and then
extending this perturbation of $\partial B^{m}$ to a perturbation of $B^{m}$ by
coming to the origin. Then $Z$ can be taken to be the final image of $X_{0}
\cap D$, which is the same as $c(X_{0} \cap\partial B^{m} \cap D_{1})$. The
arithmetic is: 
\begin{align*}
dem\, Z = dem\, c(X_{0} \cap \partial B^{m} \cap D_{1})  &  \leq1 + dem\,(X_{0} \cap
\partial B^{m} \cap D_{1})\\
&  \leq1 + dem\, (X_{0} \cap \partial B^{m}) + dem\, D_{1} - (m-1)\\
&  \leq dem\, X_{0} + dem\, A - (m-1).
\end{align*}
\indent We can easily arrange that in addition the origin $0 \notin A$, so
that $X_{0} \cap A$ lies in a truncated cone $Z_{0} = Z - int\; rB^{m}$, where
$r > 0$ is small.

\indent Let $\{ N_{\ell}= N_{\ell}(Z)\}$ be a fixed special neighborhood basis
of $Z$, (not $X_{0}!$), constructed as usual with respect to the given
coordinate chart structure, so that in particular the cone structures on
$X_{0}$ and the $N_{\ell}$'s are compatible. Let $\delta> 0$ and $p
\geq2$ be as provided by the Squeezing Proposition, for this neighborhood
sequence $\{ N_{\ell}\}$ and the given $\epsilon > 0$ (without loss $\delta<
\epsilon$, and also $(2/p)B^{m} \cap Z_{0} = \emptyset)$.

\indent The basic idea of the proof is this. For any sequence $\lambda: p
\leq\lambda(1) < \ldots< \lambda(p)$, the squeezing homomorphism $h_{\lambda}%
$, defined earlier, has the property that $h_{\lambda}(X_{0}) \subset X_{0} -
Z_{0}$, hence $h_{\lambda}(X_{0}) \cap A = \emptyset$. So, if before applying
such an $h_{\lambda}$, we can find a homeomorphism $(G_{p}$ below) of $U_{0} -
Y_{0}$ under which all $Y_{j}$-images which intersect $N_{\lambda(1)}$ are
$\delta$-small and each intersects at most one of $fr\,N_{\lambda(1)}, \ldots
fr\,N_{\lambda(p)}$, then the homeomorphism $h = h_{\lambda}G_{p}$ will satisfy
the $a$-Shrinking Lemma.

\indent The details follow. From this point on, the proof is very similar to
the preceding one.

\indent To start, let $\mathscr{Y}_{1}$ be the finite subcollection of members of $\mathscr{Y} -
\{Y_{0}\}$ which intersect $Z$ and have diameter $\geq\delta$. Choose a
collection $\mathscr{U}_{1}$ of disjoint open saturated neighborhoods of the members of
$\mathscr{Y}_{1}$, each member of $\mathscr{U}_{1}$ having diameter $< \epsilon$ and lying in $U_{0} - Y_{0}$. For each member $U$ of $\mathscr{U}_{1}$, apply the $(a-1)$-Shrinking Lemma,
with $\epsilon$-value $\min\{\delta, dist \; (\cup \mathscr{Y}_{1}, M - \cup \mathscr{U}_{1})\}$, to
find a homeomorphism $h_{U}$, supported in $U$, such that the $h_{U}$-image of
any member of $\mathscr{Y}$ lying in $U$ and intersecting $Z$ has diameter $< \delta$.
Letting $H_{1}$ be the composition of these $h_{U}$'s, $U \in \mathscr{U}_{1}$, it
follows that each member of $H_{1}(\mathscr{Y} - \{Y_{0}\})$ which intersects $Z$ has
diameter $< \delta$. So we can choose $\lambda(1) \geq p$ so large that
$N_{\lambda(1)} \subset U_{0}$ and each member of $H_{1}(\mathscr{Y} - \{Y_{0}\})$ which
intersects $N_{\lambda(1)}$ has diameter $< \delta$.

\indent From now on, the repeating steps are the same, but they are a little
bit different from the just-completed first step. In general, given $k$, $2
\leq k \leq p$, suppose we have constructed a homeomorphism $G_{k-1}(= H_{k-1}
\circ\ldots\circ H_{1}): M \rightarrow M$, supported in $U_{0} - Y_{0}$, and a
sequence $p \leq\lambda(1) < \ldots< \lambda(k-1)$ with the properties:
\\[6pt]
\indent $(1_{k-1})$ each member of $G_{k-1}(\mathscr{Y} - \{Y_{0}\})$ lying in $U_{0}$
has diameter $<\epsilon$, and each member which intersects $N_{\lambda(1)}$ has
diameter $< \delta$, and
\\[6pt]
\indent $(2_{k-1})$ each member of $G_{k-1}(\mathscr{Y} - \{Y_{0}\})$ intersects at most
one of $fr\,N_{\lambda(1)}, \ldots,\newline fr\,N_{\lambda(k-1)}$.
\\[6pt]
We show how to construct the analogous $G_{k}$ and $\lambda(k)$. Let $\mathscr{Y}_{k}$
be the finite subcollection of members of $G_{k-1}(\mathscr{Y} - \{Y_{0}\})$ which
intersect both $fr\,N_{\lambda(k-1)}$ and $Z$. Choose a collection $\mathscr{U}_{k}$ of
disjoint open $G_{k-1}(\mathscr{Y})$-saturated neighborhoods of the members of $\mathscr{Y}_{k}$,
each having diameter $< \delta \:( < \epsilon)$ and lying in $N_{\lambda(k-2)} -
Y_{0}$. (Let $N_{\lambda(0)}$ be $U_{0}$ here.) For each member $U$ of $\mathscr{U}_{k}$,
apply the $(a-1)$-Shrinking Lemma, with $\epsilon$-value $\min\{ dist(Z, M -
N_{\lambda(k-1)}), dist(\cup \mathscr{Y}_{k}, M - \cup \mathscr{U}_{k})\}$, to find a
homeomorphism $h_{U}$, supported in $U$, such that the $h_{U}$-image of any
member of $G_{k-1}(\mathscr{Y} - \{Y_{0}\})$ which intersects $Z$ necessarily misses $fr\,
N_{\lambda(k-1)}$. Letting $H_{k}$ be the composition of these $h_{U}$'s, $U
\in \mathscr{U}_{k}$, and letting $G_{k} = H_{k} \circ G_{k-1}$, it follows that $G_{k}$
satisfies properties $(1_{k})$ and $(2^{\prime}_{k})$, where $(2^{\prime}%
_{k})$ is property $(2_{k})$ with $Z$ in place of $fr N_{\lambda(k)}$. To
achieve $(2_{K})$, simply choose $\lambda(k) > \lambda(k-1)$ so large that
each member of $G_{k}(\mathscr{Y} - \{Y_{0}\})$ intersects at most one of $fr\,N_{\lambda(k-1)}$ and $fr\,N_{\lambda(k)}$.

\indent After constructing $G_{p}$ in this manner, with properties $(1_{p})$
and $(2_{p})$, the final homeomorphism $h$ of the $a$-Shrinking Lemma is $h =
h_{\lambda}G_{p}$, as explained earlier. It follows from the Squeezing
Proposition that this $h$, which is supported in $U = U_{0}$, has the desired
properties. This completes the proof that $(a - 1)$-Shrinking Lemma
$\Longrightarrow$ $a$-Shrinking Lemma.

\section{Shrinking tame closed-codimension 3 decompositions.}

This section may be regarded as a (somewhat optional) warmup for \S \ 4. The
goal here is to prove the 1-LCC Shrinking Theorem (so named by J. Cannon in
\cite{Ca}) below, using the 0-Dimensional Shrinking Theorem of \S 2. The proof
introduces the key idea of \S 4, without some of the surrounding
complications. But in as much as \S 4 uses only the 2-dimensional case of the
1-LCC Shrinking Theorem, which has been proved by Tinsley \cite{Ti} for ambient
dimension $\geq6$, the anxious reader may skip directly to \S 4.
\\\\
\textbf{1-LCC Shrinking Theorem.} Suppose $f : M \rightarrow Q$ is a cell-like
map of a manifold $M$ onto a quotient space $Q$, such that the closure in $Q$
of the image of the nondegeneracy set of has dimension $\leq m - 3$, and is
1-LCC in $Q$. Suppose $\dim M \geq5$. Then $f$ is arbitrarily closely
approximable by homeomorphisms, i.e., the decomposition of $M$ induced by $f$
is shrinkable.

\section{Proof of the Approximation Theorem.}

The basic input into this section is the 0-Dimensional Shrinking Theorem of
\S 2, and the 1-LCC Shrinking Theorem of \S 3 for the case where the closure
of the image of the nondegeneracy set is 2-dimensional (and the ambient
dimension is $\geq5$).

Let $f : M \rightarrow Q$ be as in the statement of the Approximation
Theorem. The first task is to filter $Q$ by a sequence of $\sigma$-compact
subsets, over which $f$ will be made a homeomorphism, in order of their
increasing dimension. We write $Q = P^{q} \supset P^{q-1} \supset\ldots\supset
P^{2}$ where:\\[6pt]
\indent(1) each $P^{i}$ is a $\sigma$-compact subset of $Q$, with $\dim P^{i} \leq i$
and $\dim(P^{i} - P^{i-1}) \leq0$ (hence $\dim(Q-P^{i}) \leq q - i - 1$, by \cite{Hu-Wa});\\[6pt]
\indent(2) $P^{q-3}$ is 1-LCC in $Q$, and \\[6pt]
\indent(3) any $\sigma$-compact subset of $Q-P^{2}$ is 1-LCC in $Q$.
\\\\
\indent These properties can be achieved as follows:\\[6pt]
\textit{Property (1).} One starts with $P^{q} = Q$, where $q = \dim Q = \dim M$ \cite{Ko},
and works down. (Actually, it is necessary only to use in what follows the
fact that $5 \leq q < \infty$, the former inequality to ensure that $P^{2}
\subset P^{q-3})$. Having defined a $\sigma$-compactum $P^{i}$ in $Q$ with
$\dim P^{i} \leq i$, one can let $P^{i-1}$ be the union of the frontiers (in
$P^{i}$) of a countable topology basis of open subsets of $P^{i}$, each with
frontier of dimension $\leq i - 1$, \cite{Hu-Wa}.\\[6pt]
\indent\textit{Property (2).} Let $\mathscr{A}$ be a countable dense subset of $\mbox{Maps}(B^{2}, Q)$,
the set of maps of the 2-cell $B^{2}$ to $Q$ with the uniform topology. (Recall
$\mbox{Maps}(X,Y)$, for $X,Y$ compact metric, is a complete separable metric
space.) Because $Q$ has the disjoint disc property, each map in $\mathscr{A}$ can be
chosen to have image of dimension $\leq2$, because each map in $\mathscr{A}$ can be
chosen to be an embedding. (\textit{Question}: Can these 2-cell images be chosen
2-dimensional, merely assuming $Q$ is an arbitrary compact metric $ANR$? Or
more strongly, assuming $Q$ is an $ANR$ homology manifold, or even an $ANR$
cell-like image of a manifold, but not assuming the disjoint disc property?)
Let $A$ denote the 2-dimensional union of the images of the maps in $\mathscr{A}$. To
achieve property (2), it suffices to construct $P^{q-3}$ so that $P^{q-3} \cap
A = \emptyset$. This can be done by what amounts to a relative version of the
construction for property (1). Namely, given any $\sigma$-compact 2-dimensional
subset $A$ of $Q$, one constructs $P^{q-1}$ as above so that in addition
$\dim(P^{q-1} \cap A) \leq1$, then in $P^{q-1}$ one constructs $P^{q-2}$ as
above so that in addition $\dim(P^{q-2} \cap A) \leq0$, etc. The general
dimension theory fact, from which the desired countable topology bases of open
sets can be constructed in the successive $P^{i}$'s, is the following.\\\\
\newline\textbf{Proposition:} \textit{Given any $\sigma$-compact subset A of a $\sigma$-compact
metric space $Q$, and given any point $x \in A$, then $x$ has arbitrarily
small neighborhoods $\{ U \}$ such that $\dim fr\,U \leq\dim Q - 1$ and $\dim
fr_{A}(U \cap A) \leq\dim A - 1$.}
\\\\
\indent \textbf{Note:} I would guess that the Proposition holds for an arbitrary subset
$A$ of an arbitrary separable metric space $Q$, i.e., that both occurrences of
``$\sigma$-compact'' can be dropped ($Q$ separable). But the above version is
all that is required here.
\\\\
\textit{\textbf{Proof}} (heavy-handed;perhaps it will be improved). 
Suppose $Q$ is finite-dimensional. Embed $Q$
tamely in some large dimensional euclidean space $\mathbb{R}^{n}$. The goal is, by
ambient isotopy of $\mathbb{R}^{n}$, to move $Q$ so that $A \subset N^{a}_{n}$ (= the
$a$-dimensional N\"{o}beling space in $\mathbb{R}^{n}$. The definition of $N^{\ell}_{n}$ is recalled in \S 2 in the first paragraph of the proof of the Lemma.  Or see \cite{Hu-Wa}.) and $Q \subset N^{q}_{n}$. This will establish the Proposition, for pairs $(Q,A)$ in $(N^{q}_{n}, N^{a}_{n})$
have the desired property because $(N^{q}_{n}, N^{a}_{n})$ does, as can be
verified directly. (Take cube neighborhoods with rational faces.)

\indent To move $(Q,A)$ into $(N^{q}_{n}, N^{a}_{n})$, one moves the $\sigma
$-compact complements of $N^{q}_{n}$ and $N^{a}_{n}$ off of $Q$ and $A$,
respectively. This is done as usual via a limit argument, moving ever larger
compacta of $\mathbb{R}^{n} - N^{q}_{n}$ (and of $\mathbb{R}^{n} - N^{a}_{n})$ off of ever
larger compacta of $Q$ (and of $A$). $\blacksquare$ \\\\
\indent\textit{Property (3).} At the same time one constructs the set $\mathscr{A}$
for property (2), one can construct a set $\mathcal{B}$ with precisely the
same properties as $\mathcal{A}$ (i.e., $\mathcal{B}$ is a countable dense
subset of maps$(B^{2},Q)$, with images having dimension $\leq2$), such that in
addition the set $B$ of $\mathcal{B}$-images is disjoint from the set $A$, of
$\mathcal{A}$-images. This is done in \cite{Ca}. Now, having gotten such an $A$
and $B$, and having constructed the $P^{i}$'s in the manner already
described to satisfy properties (1) and (2), one makes the $P^{i}$'s satisfy
in addition property (3), by replacing each $P^{i}$ by $P^{i} \cup B$.

\indent We can now proceed to the proof itself, which is broken into three
steps. The proof bears a curious resemblance to dual skeleton arguments used
in engulfing. The 0-dimensional Shrinking Theorem of \S 2 can be thought of as
the analogue of codimension $\geq3$ engulfing. Step II below, which I regard
as the key idea of this section, can be thought of as the analogue of the step
in dual skeleton engulfing arguments where one pushes across from the
codimension 3 skeleton toward the dual 2-skeleton.\newline\newline\textbf{Step
I.} Given $f : M \rightarrow Q$ as in the statement of the Approximation
Theorem, and given the filtration $Q = P^{q} \supset P^{q-1} \supset
\ldots\supset P^{2}$ as constructed above, the goal of this step is to
construct a cell-like map $f_{I}: M \rightarrow Q$, arbitrarily close to $f$,
such that $f_{I}$ is $1-1$ over $P^{2}$, that is, the nondegeneracy set of
$f_{I}$ misses $f^{-1}_{I}(P^{2})$. (This happens to make $f^{-1}_{I}|P^{2}$ an
embedding, but that is not of direct relevance.)

\indent To achieve Step I, we would like to say ``use the 1-LCC Shrinking
Theorem of \S 3 to shrink the decomposition of $M$ induced by the restriction
of $f$ over $P^{2}$''. However, this makes no sense, as this decomposition of
$M$ may not be uppersemicontinuous. If $P^{2}$ were compact, this would work
nicely. What we can do is to shrink this $P^{2}$-induced decomposition over
larger and larger compact subsets of the $\sigma$-compactum $P^{2}$, so that
in the limit the decomposition over $P^{2}$ is shrunk, but the decomposition
over $Q-P^{2}$ may be changed (e.g., some points of $M$ may get ``blown up''
to be nondegenerate elements). The details follow, cast in the language of
cell-like maps (as opposed to decompositions).

\indent Write $P^{2} = \cup^{\infty}_{j=1} P^{2}_{j}$, where $P^{2}_{j}$ is
compact and $P^{2}_{j} \subset P^{2}_{j+1}$. The desired map $f_{I} : M
\rightarrow Q$ of Step I is gotten by taking the limit of a sequence of
cell-like maps $\{ f_{j} : M \rightarrow Q|j \geq0\}$ which are constructed to
have the following properties $(j \geq1; f_{0} \equiv f)$:\\[6pt]
\indent $(i)$ $f_{j}$ is $\epsilon/2^{j}$-close to $f_{j-1}$, where $\epsilon>
0$ is the desired degree of closeness of $f_{I}$ to $f$,\\[6pt]
\indent $(ii)$ $f_{j}$ is 1-1 over $P^{2}_{j}$ (that is, the nondegeneracy set
of $f_{j}$ misses $f^{-1}_{j}(P^{2}_{j}))$, and\\[6pt]
\indent $(iii)$ $f_{j}$ agrees with $f_{j-1}$ over $P^{2}_{j-1}$, and $f_{j}$ is
majorant-closed to $f_{j-1}$ over $P - P^{2}_{j-1}$. In precise terms,
\begin{align*}
{dist}(f_{j}(x), f_{j-1}(x)) \leq\epsilon_{j}(x) \equiv(1/3^{j})\;{dist}\; (f_{j-1}(x), P^{2}_{j-1})
\end{align*}
for each $x \in M$. (Disregard $(iii)$ when $j = 1$.)\\[6pt]
\indent Condition $(i)$ guarantees that the $f_{j}$'s converge to a map $f_{I} :
M \rightarrow Q$; it is cell-like, being the limit of cell-like maps (c.f.
Introduction). Conditions $(ii)$ and $(iii)$ guarantee that $f_{I}$ is 1-1 over
$P^{2}$, as can easily be checked.

\indent To construct $f_{1}$, one simply ``shrinks the decomposition of $M$
induced by $f$ over $P^{2}_{1}$'', by applying \S 3. That is, let $M_{1}
\equiv M / \{f^{-1}_{0}(y)|y \epsilon P^{2}_{1} \}$ be the quotient space of
$M$ gotten by identifying to points the point-inverses of $P^{2}_{1}$ under
$f_{0}$. Then $P^{2}_{1}$ is 1-LCC in $M_{1}$, because $P^{2}_{1}$ is 1-LCC in
$Q$ (being a subset of $P^{q-3}$). By the 1-LCC Shrinking Theorem in \S 3, the
cell-like projection map $\pi_{1} : M \rightarrow M_{1}$ is arbitrarily
closely approximable by homeomorphism, $h_{1}$ say. Let $f_{1} = f_{0}
\pi^{-1}_{1} h_{1} : M \rightarrow Q$, which closely approximates $f_{0}$
because $h_{1}$ closely approximates $\pi_{1}$.

\indent To construct $f_{2} : M \rightarrow Q$, one can throw away $P^{2}_{1}$
from $Q$ and $f^{-1}_{1}(P^{2}_{1})$ from $M$, restricting ones attention to
the cell-like map $f_{1}|: M - f^{-1}_{1}(P^{2}_{1}) \rightarrow Q - P^{2}%
_{1}$. Arguing just as in the paragraph above, applying \S 3 now to the
(noncompact) manifold $M - f^{-1}_{1}(P^{2}_{1})$ and its cell-like quotient
\[ (M - f^{-1}_{1}(P^{2}_{1}))/\{f^{-1}_{1}(y)|y \in P^{2}_{2} - P^{2}_{1}\}, \]
one can construct, arbitrarily (majorant) close to $f_{1}|$, a cell-like map
$f^{\prime}_{2} : M - f^{-1}_{1}(P^{2}_{1}) \rightarrow Q - P^{2}_{1}$ which
is 1-1 over $P^{2}_{2} - P^{2}_{1}$. If the approximation is close enough, then
$f_{2} \equiv f^{\prime}_{2} \cup f_{1}|f^{-1}_{1}(P^{2}_{1}) : M \rightarrow
Q$ satisfies the desired properties. One continues this way to construct the
remaining $f_{j}$'s, hence $f_{I}$, completing Step I.
\\\\
\textbf{Step II.} The cell-like map $f_{I} : M \rightarrow Q$ constructed in
Step I, which is 1-1 over $P^{2}$, may nevertheless have nondegeneracy set (in
$M$) of large demension, even demension $m$. The goal in this step is to
arbitrarily closely approximate $f_{I}$ by a cell-like map $f_{II}: M
\rightarrow Q$ such that $f_{II}$ is 1-1 over $P^{2}$ \underline{and} the
nondegeneracy set of $f_{II}$ has codemension $\geq3$. This latter property will be achieved by making the nondegeneracy set
of $f_{II}$ lie in $M - L^{2}$, where $L^{2}$ is a certain countable union of
locally flat 2-planes in $M$ (so that $M - L^{2}$ is an analogue in $M$ of the
codimension 3 N\"{o}beling subspace of $\mathbb{R}^{m}$).

\indent To be precise, let $L^{2}(\mathbb{R}^{m}) \subset \mathbb{R}^{m}$ be the set of all
points in $\mathbb{R}^{m}$ having at least $m-2$ coordinates rational. Then
$L^{2}(\mathbb{R}^{m})$ is a countable union of 2-dimensional hyperplanes in $\mathbb{R}^{m}$,
each hyperplane being a translate of one of the $m(m-1)/2$ standard
2-dimensional coordinate subspaces of $\mathbb{R}^{m}$. N\"{o}beling's space is $\mathbb{R}^{m}
- L^{2}(\mathbb{R}^{m})$. In $M$, let $\{ \phi_{j} : \mathbb{R}^{m} \rightarrow M\}$ be a
locally finite cover by coordinate charts, and define $L^{2} = \cup_{j}
\phi_{j}(L^{2}(\mathbb{R}^{m}))$. Then, just as for N\"{o}beling's space, any compact
subset of $M - L^{2}$ has codemension $\geq3$ in $M$.

\indent Write $L^{2} = \cup^{\infty}_{j=1} L^{2}_{j}$, where each $L^{2}_{j}$
is a finite 2-complex and $L^{2}_{j} \subset L^{2}_{j+1}$. ($L^{2}_{j}$ need
not be a subcomplex of $L^{2}_{j+1}$ for the argument below.) The desired map
$f_{II} : M \rightarrow Q$ of Step II is gotten by taking the limit of a
sequence of cell-like maps $\{ f_{j} : M \rightarrow Q|j \geq 0\}$, where $f_{0} =
f_{I}$ and each $f_{j}$ is gotten from $f_{j-1}$ by preceding $f_{j-1}$ by a
homeomorphism of $M$ which moves $L^{2}_{j}$ off of the nondegeneracy set of
$f_{j-1}$. Thus each $f_{j}$ will have its nondegeneracy set qualitatively the
same as that of $f_{I}$. But as $j$ increases, the nondegeneracy set will be
getting better and better controlled by being moved off larger and larger
compact pieces of $L^{2}$, so that in the limit the nondegeneracy set
completely misses $L^{2}$ (and thus its quality may change severely, but at
least its codemension becomes $\geq3$). The other desired property of the
limit map $f_{II}$, that it remain 1-1 over $P^{2}$, will be achieved as in
Step I, by keeping the $f_{j}$'s controlled over the larger and larger compact
subsets $\{P^{2}_{j}\}$ of $P^{2}$.

\indent The precise properties of the $f_{j}$'s are $(j \geq1)$:\\[6pt]
$(i)$ $f_{j}$ is $\epsilon/2^{j}$ close to $f_{j-1}$, where $\epsilon>
0$ is the desired degree of closeness of $f_{II}$ to $f_{I} = f_{0}$\\[6pt]
$(ii)$ $f_{j}$ is 1-1 over $f_{j}(L^{2}_{j}) \cup P^{2}$ (that
is, the nondegeneracy set of $f_{j}$ misses $L^{2}_{j} \cup f^{-1}_{j} \cup
f^{-1}_{j}(P^{2}))$, and\\[6pt]
$(iii)$ $f_{j}$ agrees with $f_{j-1}$
over $f_{j-1}(L^{2}_{j-1}) \cup P^{2}_{j-1} \equiv W_{j-1}$, and $f_{j}$ is
majorant-close to $f_{j-1}$ over $Q - W_{j-1}$. In precise terms,
\[
\mbox{$dist(f_{j}(x), f_{j-1}(x)) \leq\epsilon_{j}(x) \equiv(1/3^{j})dist(f_{j-1}(x), W_{j-1})$}
\]
for each $x \in M$.  (Disregard $(iii)$  when j = 1.)\\[6pt]
\indent As in Step I, the reader can verify that these properties ensure that
the limit map $f_{II}: M \rightarrow Q$ has the desired properties. So it
remains to explain how these properties of the $f_{j}$'s are achieved.

\indent Consider $f_{1}$. To construct it, we find a homeomorphism $h_{1} : M
\rightarrow M$ such that $f_{0} h_{1}$ is $\epsilon/2$ close to $f_{0}$, and
such that $h_{1}(L^{2}_{1}) \cap$ nondegeneracy set $(f_{0}) = \varnothing$.
Then we can define $f_{1} = f_{0} h_{1}$. To find $h_{1}$, the key is first to
find the tame ($\equiv\;$ 1-LCC) embedding $h_{1}|L^2_1$, call it $\alpha_{1} :
L^{2}_{1} \rightarrow M$, such that $f_{0} \alpha_{1} : L^{2}_{1} \rightarrow
Q$ is close to $f_{0}|:L^{2}_{1} \rightarrow Q$ and such that $\alpha
_{1}(L^{2}_{1})$ misses the nondegeneracy set of $f_{0}$. This embedding
$\alpha_{1}$ will be gotten by working in the quotient space $Q$. The point
is, the image in $Q$ of the nondegeneracy set of $f_{0}$ misses $P^{2}$, hence
is a 1-LCC $\sigma$-compactum, and hence $f_{0}|:L^{2}_{1} \rightarrow Q$ can
be approximated arbitrarily closely by a 1-LCC embedding $\beta_{1} :
L^{2}_{1} \rightarrow Q$ whose image misses $f_{0}(nondeg(f_{0}))$. This is a
standard argument (c.f. Introduction). Let $\alpha_{1} = f^{-1}_{0} \beta_{1} : L^{2}_{1} \rightarrow M$,
which is a 1-LCC embedding. By the usual cell-like map arguments, the small
homotopy joining $f_{0}|L^{2}_{1}$ and $\beta_{1} : L^{2}_{1} \rightarrow Q$
in $Q$ can be lifted, as efficiently as desired (efficiency being measured by
smallness in $Q$), to a homotopy joining the inclusion $L^{2}_{1}
\hookrightarrow M$ and the embedding $\alpha_{1} : L^{2}_{1} \rightarrow M$.
If $\dim M \geq 6$, this homotopy can be covered as efficiently as desired by
an ambient isotopy of $M$, whose end homomorphism $h_{1}$ restricts on
$L^{2}_{1}$ to $\alpha_{1}$. 

\indent If $\dim M = 5$, some additional remarks are called for.
\[ 
\vdots 
\]
\indent In general, to construct $f_{j}$ given $f_{j-1}, j \geq2$, one uses
the same simple device as in Step I, namely temporarily throwing away
$f_{j-1}(L^{2}_{j-1}) \cup P^{2}_{j-1} \equiv W_{j-1}$ (see property $(iii)$)
and its preimage under $f_{j-1}$. That is, one considers the cell-like
restriction map $f_{j-1}|:M - f^{-1}_{j-1}(W_{j-1}) \rightarrow Q - W_{j-1},$
and one constructs a homeomorphism $h^{\prime}_{j}$ of the source manifold $M
- f^{-1}_{j-1}(W_{j - 1})$ onto itself such that $f^{\prime}_{j-1} \equiv
f_{j-1}h^{\prime}_{j}$ is arbitrarily close to $f_{j-1}|$, and $h^{\prime}%
_{j}(L^{2}_{j} - L^{2}_{j-1}) \cap nondeg(f_{j-1}|) = \varnothing$. This of
course requires the noncompact, majorant-controlled versions of the arguments
used in the preceding paragraphs to construct $h_{1}$, but they all are
available. If the approximation of $f^{\prime}_{j-1}$ to $f_{j-1}|$ is
sufficiently close, then \vspace{6pt}
\[ 
f_{j} \equiv f^{\prime}_{j} \cup f_{j-1}|f^{-1}_{j-1}(W_{j-1}) : M \rightarrow Q 
\]
satisfies the desired properties. This completes Step II. \\\\
\textbf{Step III.} In this step, the decomposition of $M$ induced by
$f_{II} : M \rightarrow Q$ is shrunk over $P^{3}$, then over $P^{4}, \ldots,$
and finally over $P^{q} = Q$, at each stage using the 0-Dimensional Shrinking
Theorem of \S 2. To be a little bit more precise, we start this step with the
cell-like map $f_{II,2} \equiv f_{II} : M \rightarrow Q$ produced in Step II,
which is already $1-1$ over $P^{2}$ and has nondegeneracy set of codemension
$\geq3$, and we produce successive approximations $f_{II,i} : M \rightarrow Q$,
for $i$ running from 3 up to $q = \dim Q$, where each $f_{II,i}$ is
arbitrarily close to $f_{II,i-1}$, $f_{II,i}$ is $1-1$ over $P^{i}$ and
$f_{II,i}$ has nondegeneracy set of codemension $\geq3$.

\indent The details follow. Fix $i, 3 \leq i \leq q$. Given a cell-like map
$f_{II,i-1} : M \rightarrow Q$ which is $1-1$ over $P^{i-1}$ and has
nondegeneracy set of codemension $\geq3$, we show how to produce a
corresponding cell-like map $f_{II,i}: M \rightarrow Q$, arbitrarily close to
$f_{II,i-1}$. Write $P^{i} = \cup^{\infty}_{j=1}P^{i}_{j}$, where $P^{i}_{j}$
is compact and $P^{i}_{j}$ is compact and $P^{i}_{j} \subset P^{i}_{j+1}$. As
in the previous two steps, the map $f_{II,i}$ will be gotten as the limit of a
sequence of cell-like maps $\{ f_{j}:M \rightarrow Q| j \geq0\}$. The
properties of the $f_{j}$'s are $(j \geq1; f_{0} = f_{II, i-1}):$\\[6pt]
\indent $(i)$ $f_{j}$ is $\epsilon/2^{j}$-close to $f_{j-1}$,
where $\epsilon> 0$ is the desired degree of closeness of $f_{II,i}$ to
$f_{II,i-1} = f_{0}$,\\[6pt] 
\indent $(ii)$ $f_{j}$ is $1-1$ over
$f_{j}(L^{2}_{j}) \cup P^{i-1} \cup P^{i}_{j}$ (that is, the nondegeneracy set
of $f_{j}$ misses $L^{2}_{j} \cup f^{-1}(P^{i-1} \cup P^{i}_{j}))$, and
$f_{j}$ has nondegeneracy set of codemension $\geq3$ (Note: the reason for
using the $L^{2}_{j}$'s here is to control the codemension in the limit, as in
Step II), and\\[6pt]
\indent $(iii)$ $f_{j}$ agrees with $f_{j-1}$ over
$f_{j-1}(L^{2}_{j-1}) \cup P^{i}_{j-1} \equiv W_{j-1}$, and $f_{j}$ is
majorant-close to $f_{j-1}$ over $Q-W_{j-1}$. In precise terms, 
\[
dist(f_{j}(x), f_{j-1}(x)) \leq\epsilon_{j}(x) \equiv(1/3^{j})dist(f_{j-1}(x), W_{j-1})
\]
for each $x \in M$. (Disregard $(iii)$ when $j = 1$.)\\[6pt]
\indent As in Steps I and II, the reader can verify that these properties
ensure that the limit map $f_{II,i}: M \rightarrow Q$ has the desired properties, and so in particular $f_{II,q}$ is the desired homeomorphism of the
Theorem. It remains to explain how the properties of the $f_{j}^{\prime}s$ are achieved.

\indent To construct $f_{1}$, one ``shrinks the decomposition of $M$ induced
by $f_{0}$ over $P^{i}_{1}$'', using the 0-Dimensional Shrinking Theorem of
\S 2 and the fact that this decomposition being already trivial over $P^{i-1}$,
is therefore 0-dimensional. In more detail, let $M_{1} \equiv M/ \{f^{-1}_{0}
(y)|y \in P^{i}_{1}\}$ be the quotient space of $M$ gotten by identifying
to points the point-inverse of $P^{i}_{1}$ under $f_{0}$. This decomposition
of $M$ is in fact 0-dimensional, since the image of the nondegeneracy set
lies in $P^{i}_{1} - P^{i-1}$. Hence by \S 2, the cell-like quotient map
$\pi_{1} : M \rightarrow M_{1}$ is arbitrarily closely approximable by
homeomorphism, $h_{1}$ say. Define $f^{*}_{1} = f_{0} \pi^{-1}_{1} h_{1} : M
\rightarrow Q$, which closely approximates $f_{0}$ because $h_{1}$ closely
approximates $\pi_{1}$. The nondegeneracy set of $f^{*}_{1}$ has codemension
$\geq3$, because it equals $h^{-1}_{1} \pi_{1}(nondeg(f_{0}) - f^{-1}%
_{0}(P^{i}_{1}))$, and $h^{-1}_{1} \pi_{1}$ is an embedding on the open
neighborhood on $M - f^{-1}_{0}(P^{i}_{1})$ of $nondeg(f_{0}) - f^{-1}%
_{0}(P^{i}_{1})$, thus preserving its codemension. To complete this stage, let
$g_{1}: M \rightarrow M$ be a homeomorphism, arbitrarily close to the
identity, such that $g_{1}(L^{2}_{1}) \cap nondeg(f^{*}_{1}) = \varnothing$.
Then define $f_{1} = f^{*}_{1} g_{1}$.

\indent In general, to construct $f_{j}$ given $f_{j-1}, j \geq2$, we use the
now-familiar device of working in a restricted open subset of $M$. Letting $W_{j-1} \equiv f_{j-1}(L^{2}_{j-1}) \cup P^{i}_{j-1}$, we focus for
the moment on the restricted map $f_{j-1}| : M - f^{-1}_{j-1}(W_{j-1})
\rightarrow Q - W_{j-1}$. Let
\[
M^{\prime}_{j} \equiv(M - f^{-1}_{j-1}(W_{j-1}))/\{f^{-1}_{j-1}(y)|y \in
P^{i}_{j} - W_{j-1}\}
\]
be the quotient space of $M - f^{-1}_{j-1}(W_{j-1})$ gotten by identifying to
points the point-inverses of $P^{i}_{j} - W_{j-1}$ under $f_{j-1}|$. As above,
this decomposition is 0-dimensional, and so by \S 2 the cell-like quotient map
$\pi^{\prime}_{j} : M - f^{-1}_{j-1}(W_{j-1}) \rightarrow M^{\prime}_{j}$ is
arbitrarily (majorant) closely approximable by homeomorphism, $h^{\prime}_{j}$
say. Define $f^{*{\prime}}_{j} = f_{j-1} \pi^{\prime -1}_{j} h^{\prime}_{j} : M - f^{-1}_{j-1}(W_{j-1}) \longrightarrow Q-W_{j-1}$, which closely approximates $f_{j-1}|$. As above, $nondeg(f^{*{\prime}}_{j})$ has
codemension $\geq3$, and so there is a homeomorphism $g^{\prime}_{j}$ of $M -
f^{-1}_{j-1}(W_{j-1})$ onto itself, arbitrarily close to the identity, such
that $g^{\prime}_{j}(L^{2}_{j} - f^{-1}_{j-1}(W_{j-1})) \cap nondeg(f^{*{\prime}}_{j}) = \varnothing$. 
Define
\[
f^{\prime}_{j} = f^{*{\prime}}_{j} g^{\prime}_{j} : M - f^{-1}_{j-1}(W_{j-1})
\rightarrow Q - W_{j-1}.
\]
If the approximations were chosen small enough, so that $f^{\prime}_{j}$ is
sufficiently close to $f_{j-1}|$, then $f_{j} \equiv f^{\prime}_{j} \cup
f_{j-1}|f^{-1}_{j-1}(W_{j-1}): M \rightarrow Q$ satisfies the desired
properties. This completes Step III, and hence the proof of the Approximation
Theorem.
\\\\

\end{document}